\documentclass{article} 
\usepackage{iclr2021_conference,times}

\usepackage{amsmath,amsfonts,bm}

\def\eqref#1{equation~\ref{#1}}

\def\1{\bm{1}}

\def\mP{{\bm{P}}}

\def\mR{{\bm{R}}}

\DeclareMathAlphabet{\mathsfit}{\encodingdefault}{\sfdefault}{m}{sl}
\SetMathAlphabet{\mathsfit}{bold}{\encodingdefault}{\sfdefault}{bx}{n}

\usepackage{graphicx}
\usepackage{hyperref}
\usepackage{url}

\usepackage{caption}
\usepackage{subcaption}

\title{Multigrid solver with super-resolved interpolation}

 \author{ 
 Francisco Holguin \\
 University of Michigan Ann Arbor, Department of Astronomy \\
 Los Alamos National Laboratory
 \And
 GS Sidharth \\
 Los Alamos National Laboratory 
 \And
 Gavin Portwood \\
 Lawrence Livermore National Laboratory 
 }

\iclrfinalcopy

\begin{document}

\maketitle
\vspace{-1em}
\begin{abstract}
The multigrid algorithm is an efficient numerical method for solving a variety of elliptic partial differential equations (PDEs). The method damps errors at progressively finer grid scales, resulting in faster convergence compared to standard iterative methods such as Gauss-Seidel. The prolongation, or coarse-to-fine interpolation operator within the multigrid algorithm lends itself to a data-driven treatment with ML super resolution, commonly used to increase the resolution of images. We (i) propose the novel integration of a super resolution generative adversarial network (GAN) model with the multigrid algorithm as the prolongation operator and (ii) show that the GAN-interpolation improves the convergence properties of the multigrid in comparison to cubic spline interpolation on a class of  multiscale PDEs typically solved in physics and engineering simulations.
\end{abstract}

\section{Introduction}
\vspace{-1em}
Partial differential equations describe the behavior of many physical phenomena, from the smallest quantum scales to the largest scales in the universe. Efficient numerical methods are crucial to simulate these complex phenomena with current and foreseeable computational resources. The application of machine learning to accelerate and improve the accuracy of physical simulations is an increasingly more active area of research, leveraging the surge of advancements in machine learning architectures over the last decade \cite[e.g.][]{wei2018, peurifoy2018,ranade2020,sanchez2020}. While replacing conventional numerical solvers with machine learning models is promising, improving heuristic operators within existing  mathematically-formally derived numerical methods can allow for interpretable computational gains, easier implementation and more rapid deployment. 

\section{Proposal}
\vspace{-1em}

\subsection{Multigrid method}  %
\vspace{-.5em}

Multigrid methods are particularly effective for solving multiscale differential equations relevant to physics and engineering \citep[e.g.][]{trottenberg2001,gs2020multiscale}. Multigrid methods formally scale linearly with the number of unknowns in a discretized physical system, which make them attractive as high performance computing increasingly enables higher-fidelity and larger scale computations. Generally, the favorable scaling of multigrid methods is due to their ability to iteratively reduce solution errors at multiple scales by obtaining coarse solutions then utilizing the coarsened solutions to solve increasingly higher-fidelity solutions.  The utilization of coarse fidelity solutions to accelerate finer fidelity solutions is enabled by a prolongation, or interpolation, operator.

\begin{figure}
    \centering
    \begin{subfigure}[t]{0.65\textwidth}
        \includegraphics[width=\textwidth]{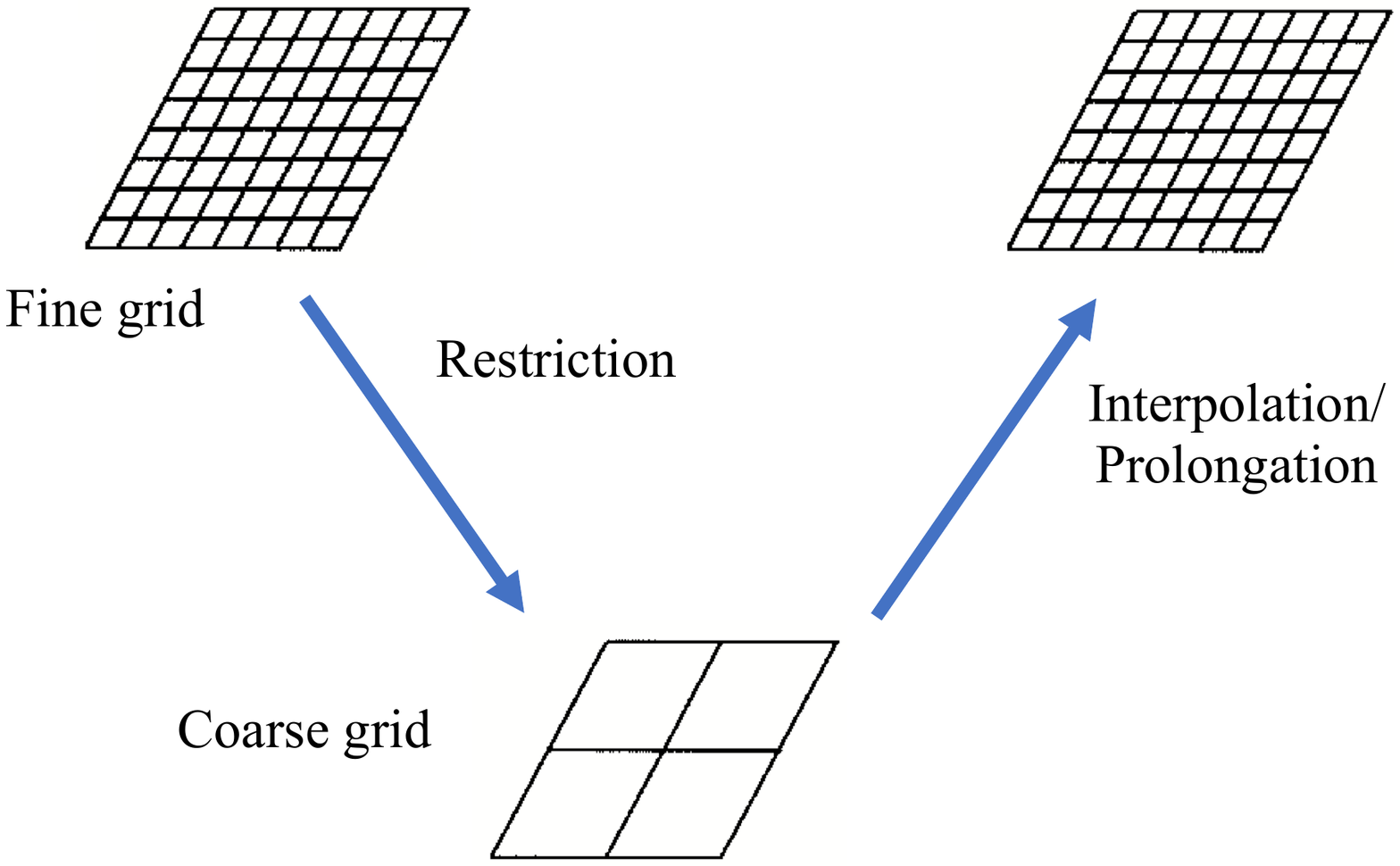}
        \caption{Schematic of a two-level multigrid cycle. Adapted from \citet{chen2001}.}
        \label{vcycle}
    \end{subfigure}
    ~
    \begin{subfigure}[t]{0.65\textwidth}
        \includegraphics[width=\textwidth]{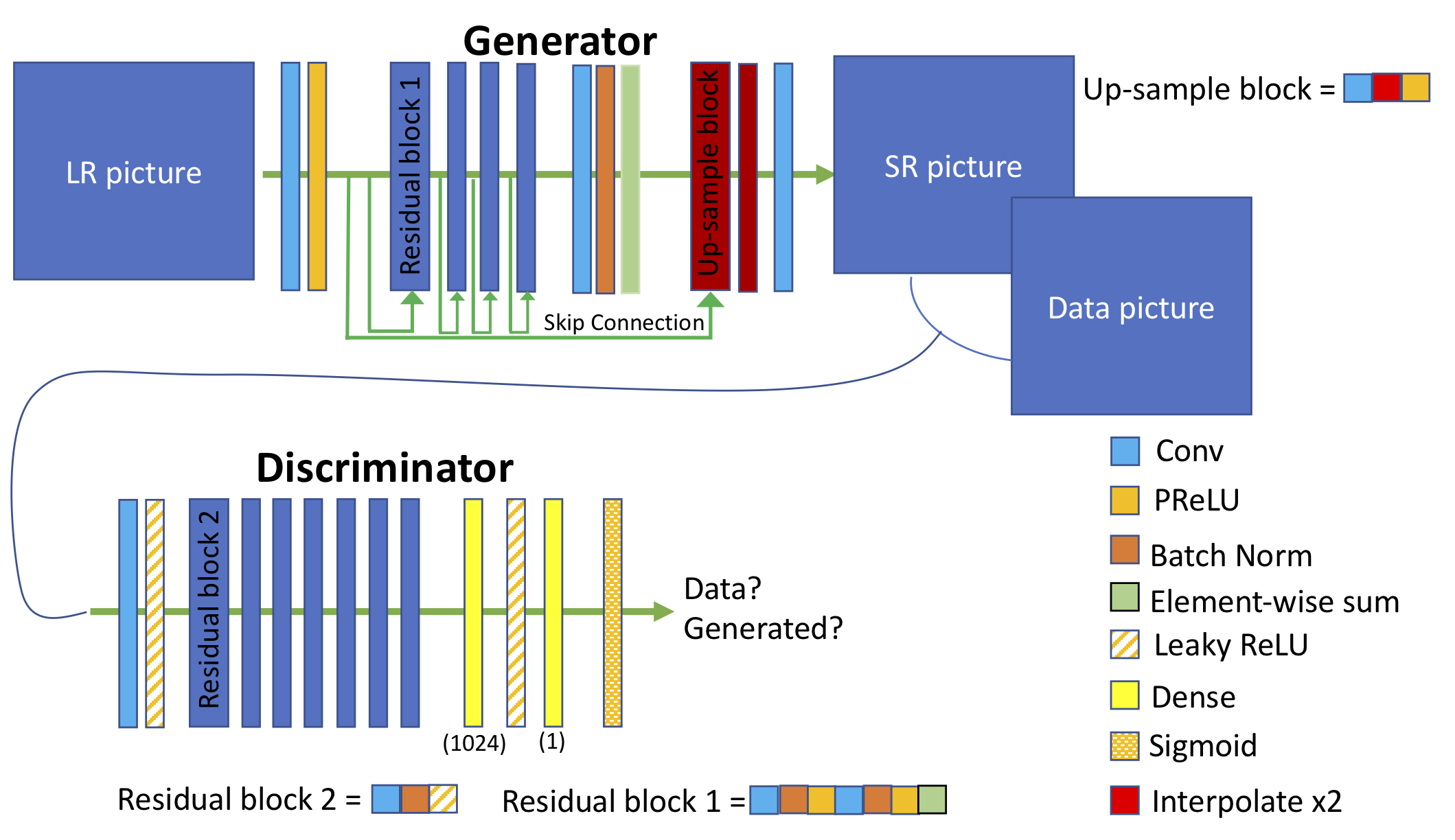}
        \caption{Structure of discriminator and generator we implement, modified from \citet{ledig2017}. }
        \label{gan_structure}
    \end{subfigure}
    \caption{Diagrams for two-level multigrid and GAN interpolation.}
\end{figure}

Prolongation within a multigrid algorithm lends itself to a data-driven treatment due to the ability for machine learning techniques to perform data-informed interpolation, as seen frequently seen in super-resolution methods as used in image analysis \citep[c.f.][]{farsiu2004advances,park2003super}. In this work, we consider the applications of enhancing multigrid methods by using a data-driven prolongation operator akin to super-resolution operators used in image analysis. In particular we consider a common Poisson differential equation, frequently solved in engineering and physics calculations, written as
\begin{equation}
    \nabla^2 p = \frac{\partial^2 p}{\partial^2 x} + \frac{\partial^2 p}{\partial^2 y} = f(x,y),
\label{poissoneq}
\end{equation}
where the solution $p$ is discretized on a grid. Figure \ref{vcycle} shows a schematic of the two-level multigrid, where the top horizontal level represents the highest resolution grid $h_2$. Applying the Poisson operator $\nabla^2 p - f$ to the grid $h_2$, using a relaxation method such as Gauss-Seidel, results in an residual $r_f$ relative to the true solution. For these relaxation methods,  high frequency errors (relative to the grid scale) are quickly reduced, while low spatial frequency errors remain. It can be shown \citep[][]{trottenberg2001} that the residual also satisfies the Poisson equation, so we recursively solve the Poisson equation for this error. In order to efficiently suppress low frequency errors, we apply a restriction operator $\displaystyle \mR$ to coarsen the grid, so that previously low spatial frequency errors on the fine grid become high frequency errors on the coarse grid. Once the coarsest level is reached, we apply a prolongation operator $\displaystyle \mP$ to the solution in order to correct the solution on the fine grid. In a uniformly spaced grid, the restriction operator $\displaystyle \mR$, can be a simple as taking only the black points in red-black ordering. The prolongation operator requires more work, as it involves an interpolation scheme to create a higher resolution grid. For example, this interpolation scheme can be a two-dimensional mesh interpolation, such as the bivariate spline function in the SciPy Python library \citep[][]{virtanen2020}. 
\subsection{GAN and integration within multigrid}
\vspace{-.5em}
\label{GAN}

Deep learning \citep{lecunhinton2015} methods have been able to capture and reproduce the structure of extremely complex physical systems \citep[e.g.][]{portwood21,wang20,he2019}. 
Even without any prior knowledge of the system, the network can learn the behavior and patterns with adequate data, relatively quickly. 
Super resolution refers to the application of these networks to low resolution images in order to predict a high resolution version. A GAN \citep{goodfellow2014} is particularly effective as an accurate super resolution model. Here, we use the GAN architecture from \citet{ledig2017}, following the example in \citet{birla2018}, as shown in Figure \ref{gan_structure}. Unlike \citet{ledig2017} who use the pre-trained VGG loss designed for photographic image recognition, we use mean squared error (MSE) loss function for the generator.


The GAN takes $n_s \times n_s$ (here $n_s=6$) point input tiles from the coarse solution
for spatially-invariant training that generates more training samples from limited datasets. Using the windowed inputs, the GAN can operate on grids with different resolution without any change in input formats. 
Ghost cells are implemented by splitting the domain up into overlapping windows and utilize only a part of the any particular GAN's predictions for the high resolution grid. 

The GAN prolongation operator works as follows:

\begin{itemize}
    \item Find the max and min data of the grid $p_{\rm{max}}, p_{\rm{min}}$, and normalize the grid data using a symmetric log function from $\displaystyle [- \infty, \infty]$ to $\displaystyle [- 1, 1]$, such that for positive data, $ \rm{log} (p_{\rm{max, min}}) = (1, 0)$ and for negative data, $ \rm{log}\   p_{\rm{max, min}} = (-1, 0)$.
    \item Divide the normalized coarse grid into overlapping $n_s \times n_s$ windows (i.e. assuming a $n_s^2$ kernel with stride = 2). Apply GAN interpolation to produce a set of high resolution $n_l^2$ window kernels. Here $n_l=24$.
    \item Construct the fine grid by assigning the central $(n_s+2)^2$ fine window values to the fine grid. Any non-overlapping areas (i.e the edges) are assigned to the window covering that area.
    \item Transform and scale back to the original pressure data.
\end{itemize}

\subsection{Data and training of the GAN}
\vspace{-.5em}
\label{training}

The pressure-Poisson formulation of the incompressible Navier-Stokes equations is a  physical system with the form of Eq. \ref{poissoneq}. The pressure-Poisson equation is used to ensure mass conservation in incompressible flows and is derived from taking the divergence of the fluid momentum equation 
\begin{equation}
    \nabla^2 p = \nabla \cdot (\nu \nabla^2 u  - ( u \cdot \nabla ) u    )  = f(x,y)
    \label{poisson_ns}
\end{equation}

Here $\nu$ is the fluid viscosity, $u$ is the fluid velocity, and $p$ is the fluid pressure. 
Evolving the system in time involves calculating the source term $f(x, y)$ from the velocity field $u$, solving the Poisson equation for pressure (e.g. using a multigrid method), and then using that pressure to update the velocity field. We produce a training set of 200 pressure fields on a grid size $192^2$ by evolving a biperiodic velocity field from a broadband initial condition using the incompressible Navier-Stokes.

We construct the training set of 1000 pressure grids as follows:
\begin{itemize}
    \item Select one of the 200 pressure fields randomly and restrict (downsample) the grid by a random power $l$ of 2 on each grid length, such that the smallest grid size is $(n_s=6)^2$. The result is a random grid of size $(192 / 2^l)^2 $.
    \item From the grid above, pick a random $6^2$ window. The total size of the space of such windows is $17,000$.
    \item Transform and rescale the window data range from $\displaystyle [- \infty, \infty]$ to $\displaystyle [- 1, 1]$ in the same manner as described in the GAN prolongation method in \ref{GAN}, using the global min and max of all the pressure grids. We used $  |p_{\rm{min}}| = 10^{-10}$ and $|p_{\rm{max}}| = 10^{-3}$.
\end{itemize}

\section{Proposal implementation}
\vspace{-1em}
\begin{figure}
    \centering
    \includegraphics[width=\textwidth]{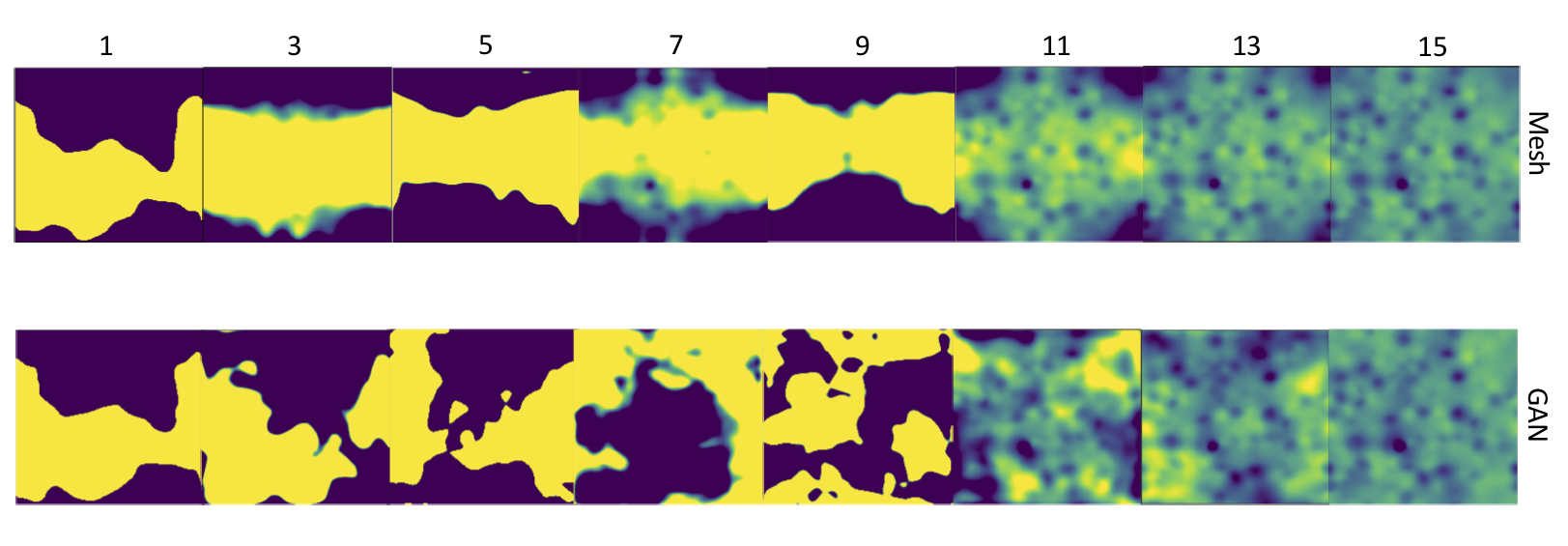}
    \caption{Comparison of the first 15 iterations of solving Eq. \ref{poissoneq} with the multigrid solver starting from an initially random grid. Top: interpolation with traditional mesh interpolation. Bottom: interpolation with GAN super resolution. }
    \label{ml_mg_data_series}
\end{figure}
\begin{figure}
    \centering
    \includegraphics[width=\textwidth]{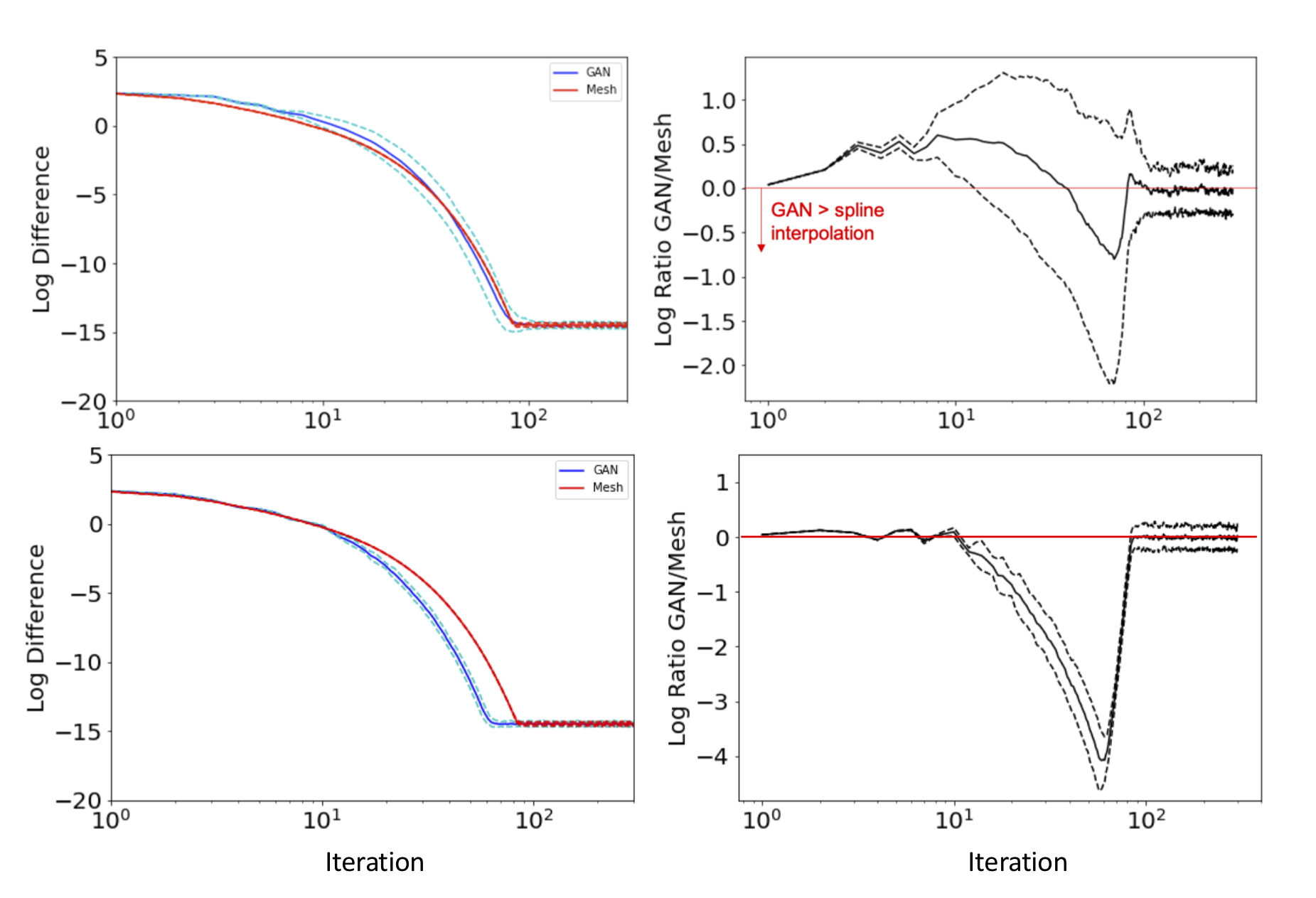}
    \caption{Averaged results of multigrid algorithm solving $100$ different grids. Norm of the difference of grid between iterations as a function of iteration, shown for two choices of prolongation operator, two-dimensional spline and GAN super resolution interpolation. The multigrid parameters used are $N_{\rm{smooth, pre}}$ = 10, $N_{\rm{smooth}} $ = 20,  $N_{\rm{step}}$ = 4, and $r_{\rm{min}}$ = 12. Top: The GAN operator is used for interpolation in every iteration. Bottom: The GAN and spline operators are alternated every other iteration.}
    \label{ml_mg_example}
\end{figure}
\begin{figure}[!ht]
    \centering
    \includegraphics[width=\textwidth]{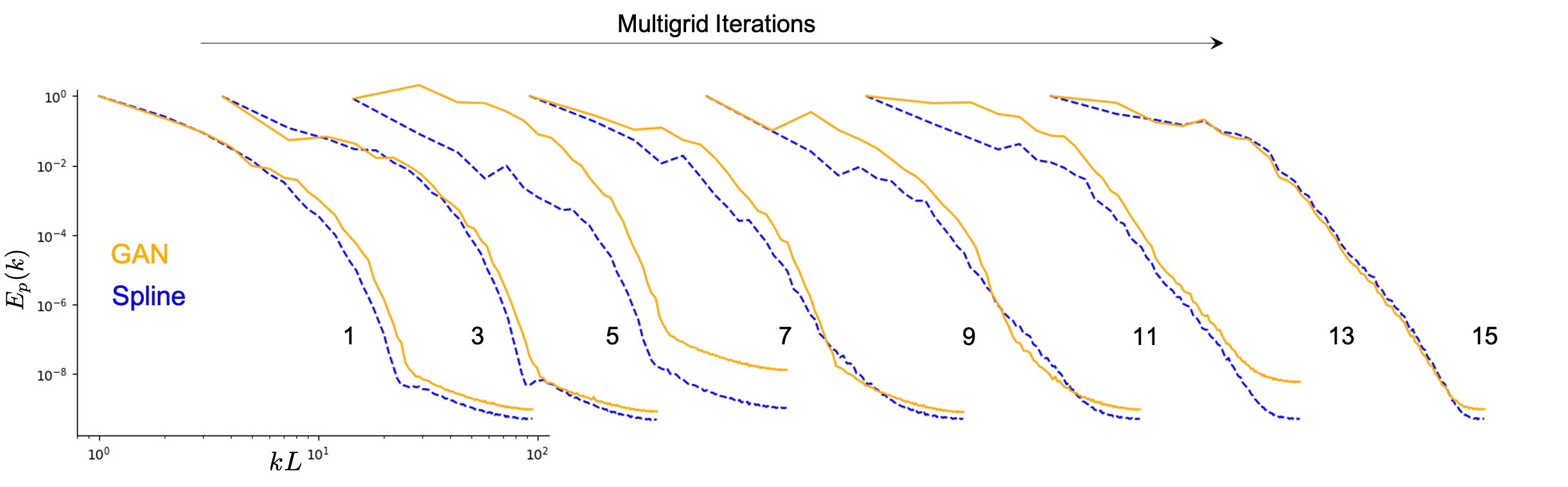}
    \caption{Convergence in $p$ spectral density with increasing multigrid iterations.}
    \label{iter}
\end{figure}
We solve Eq. \ref{poisson_ns} using a test set source grid. Our implementation of the multigrid algorithm is characterized by four parameter choices: $N_{\rm{smooth, pre}}$ = number of smoothing operations on the initial grid; $N_{\rm{smooth}}$ = number of smoothing operations at the beginning of each iteration; $N_{\rm{step}}$ = number of times the the grid resolution is reduced by a factor of two on each dimension; $r_{\rm{min}}$ = coarsest resolution. Pre-smoothing is often necessary in order to ensure convergence. 
To study the efficacy of a GAN prolongation operator under controlled settings, we implement a two-level V-cycle multigrid algorithm (fine scale $192^2$ grid to a coarse scale $12^2$ grid, corresponding to $N_{\rm{step}}$ = 5 and $r_{\rm{min}} = 12$). A two-level multigrid method with a coarsening factor of 16 is expected to perform poorly, limited by prolongation and interpolation, and therefore serves as an ideal test bed for evaluating the gains made with a super-resolution GAN. We sample $100$ pressure and source term grid pairs from our training set, following the method described in Section \ref{training}. 

Figure \ref{ml_mg_example} shows a proof-of-concept test of our multigrid solver on $100$ grids. We plot the mean profiles and $1-\sigma$ bounds  of the norm difference between the pressure grids at the current and previous iteration. The top row of the figure shows the result of using the GAN operator at every iteration compared to the spline operator. The multigrid solver converges in 100 iterations (compared to $10^4$ iterations for a simple iterative solver). Statistically, we find that using the super resolution prolongation operator can result in faster convergence compared to bivariate spline operator; the large $1-\sigma$ lines on the ratio between the two convergence lines shows that the performance is distributed. For some fields, super resolution exhibits faster convergence, but for other fields, the spline interpolation can fare better. The bottom row shows the result of alternating the GAN and spline operators compared with using only the spline operator. There is faster convergence using the alternating operators for all of the fields tested. This improvement is likely due to the fact that the GAN has been trained on real pressure field solutions and not random fields. Therefore, as soon as the pressure field evolves from random to one with physical properties, the advantages of the GAN prolongation operator set in. In practice, the initial conditions will not necessarily be random, but will be physically realistic, such as in a time-stepping physical simulation. 

These tests demonstrate strongly that there is potential for super-resolution to improve the efficiency of the multigrid algorithm by data-driven prolongation operators. Figure \ref{ml_mg_data_series} shows a sequence of the first 15 iterations of a particular field, for the solver shown in the top row of Figure \ref{ml_mg_example}. While both solvers converge to machine zero, with the conventional spline-based interpolation, the early iterations are unable to capture mid to high-frequency information. With the GAN-based interpolation however, the solver feeds in higher frequency information in the early iterations. This is confirmed via the pressure field spectra shown in Figure \ref{iter}, which show compare the spectra for the first 15 iterations.

\section{Conclusion}
\vspace{-1em}
In this research, we consider novel integration of data-driven interpolation into the robust multigrid differential equation solution technique as applied to a common differential system utilized in many physics and engineering domains. We have demonstrated the ability for super-resolution techniques, GANs in particular, to serve as a prolongation operator which can improve solution time when compared to high order spline interpolation techniques. As we continue our investigations for this approach, we propose that the integration of data-driven methods into formally-derived numerical solution procedures is a promising approach for computational acceleration of simulations.
\\

\bibliography{iclr2021_conference}
\bibliographystyle{iclr2021_conference}


\end{document}